\documentclass{amsart}
\usepackage{amssymb}
\usepackage{a4}
\begin{document}
\newtheorem{theorem}{Theorem}[section]
\newtheorem{definition}[theorem]{Definition}
\newtheorem{lemma}[theorem]{Lemma}
\newtheorem{example}[theorem]{Example}
\newtheorem{remark}[theorem]{Remark}
\newtheorem{corollary}[theorem]{Corollary}
\def\ffrac#1#2{{\textstyle\frac{#1}{#2}}}
\def\qedbox{\hbox{$\rlap{$\sqcap$}\sqcup$}}
\font\pbglie=eufm10\def\gg{\text{\pbglie g}}\def\oo{\text{\pbglie o}}
\def\id{\text{Id}}
\makeatletter
 \renewcommand{\theequation}{%
 \thesection.\alph{equation}}
 \@addtoreset{equation}{section}
 \makeatother
\title[Generalized plane wave manifolds]
{Generalized plane wave manifolds}
\author{P. Gilkey and S. Nik\v cevi\'c}
\begin{address}{PG: Mathematics Department, University of Oregon,
Eugene Or 97403 USA.\newline Email: {gilkey@darkwing.uoregon.edu}}
\end{address}
\begin{address}{SN: Mathematical Institute, SANU,
Knez Mihailova 35, p.p. 367,
11001 Belgrade,
Serbia and Montenegro.
\newline Email: {\it stanan@mi.sanu.ac.yu}}\end{address}
\begin{abstract}
   We show that generalized plane wave manifolds are complete, strongly geodesically convex, Osserman,
   Szab\'o, and Ivanov-Petrova. We show their holonomy groups are nilpotent and that all the local Weyl scalar invariants of these
manifolds vanish. We construct
   isometry invariants on certain families of these manifolds which are not of Weyl type. Given $k$, we exhibit
   manifolds of this type which are $k$-curvature homogeneous but not locally homogeneous. We also construct a
manifold which is weakly $1$-curvature homogeneous but not $1$-curvature homogeneous.
\end{abstract}
\maketitle
\section{Introduction} We begin by introducing some notational conventions. Let $\mathcal{M}:=(M,g)$ where
$g$ is a pseudo-Riemannian metric of signature $(p,q)$ on smooth manifold $M$ of dimension $m:=p+q$.
\subsection{Geodesics} We say
that
$\mathcal{M}$ is {\it complete} if all geodesics extend for infinite time and that $\mathcal{M}$ is
{\it strongly geodesically convex} if there exists a unique geodesic
between any two points of $M$;  if $\mathcal{M}$ is complete and strongly geodesically convex, then the exponential map is a
diffeomorphism from $T_PM$ to $M$ for any $P\in M$.

\subsection{Scalar Weyl invariants} Let $\nabla^kR$ be the $k^{\operatorname{th}}$ covariant derivative of the curvature operator
defined by the Levi-Civita connection.  Let $x:=(x_1,...,x_m)$ be local coordinates on $M$. Expand
\begin{equation}\label{eqn-1.a}
\nabla_{\partial_{x_{j_1}}...}\nabla_{\partial_{x_{j_l}}}R(\partial_{x_{i_1}},\partial_{x_{i_2}})\partial_{x_{i_3}}
    =R_{i_1i_2i_3}{}^{i_4}{}_{;j_1...j_l}\partial_{x_{i_4}}
\end{equation}
where we adopt the Einstein convention and sum over repeated indices. Scalar invariants of the
metric can be formed by using the metric tensors
$g^{ij}$ and $g_{ij}$ to fully contract all indices.  For
example, the scalar curvature $\tau$, the norm of the Ricci tensor $|\rho|^2$, and the norm of the full curvature
tensor $|R|^2$ are given by
\begin{eqnarray}
&&\tau:=g^{ij}R_{kij}{}^k,\nonumber\\
&&|\rho|^2:=g^{i_1j_1}g^{i_2j_2}R_{ki_1j_1}{}^kR_{li_2j_2}{}^l,\quad\text{and}\label{eqn-1.b}\\
&&|R|^2:=g^{i_1j_1}g^{i_2j_2}g^{i_3j_3}g_{i_4j_4}R_{i_1i_2i_3}{}^{i_4}R_{j_1j_2j_3}{}^{j_4}\,.\nonumber
\end{eqnarray}
Such invariants are called {\it Weyl invariants}; if all possible such invariants vanish, then $\mathcal{M}$
is said to be {\it VSI} (vanishing scalar invariants). We refer to Pravda, Pravdov\'a, Coley, and
Milson \cite{PPCM02} for a further discussion.

\subsection{Natural operators defined by the curvature tensor} If $\xi$ is a tangent vector, then
the {\it Jacobi operator}
$J(\xi)$ and the {\it Szab\'o operator}
$\mathcal{S}(\xi)$ are the self-adjoint linear maps which are defined by:
$$
J(\xi):x\rightarrow R(x,\xi)\xi\quad\text{and}\quad\mathcal{S}(\xi):x\rightarrow\nabla_\xi R(x,\xi)\xi\,.
$$
Similarly if
$\{e_1,e_2\}$ is an oriented orthonormal basis for an oriented spacelike (resp. timelike) $2$-plane $\pi$,
the {\it skew-symmetric curvature operator} $\mathcal{R}(\pi)$ is defined by:
$$\mathcal{R}(\pi):x\rightarrow R(e_1,e_2)x\,.$$

\subsection{Osserman, Ivanov-Petrova, and Szab\'o manifolds} We say that $\mathcal{M}$ is {\it spacelike Osserman} (resp. {\it
timelike Osserman}) if the eigenvalues of
$J$ are constant on the pseudo-sphere bundles of unit spacelike (resp. timelike) tangent vectors. The notions
{\it spacelike Szab\'o}, {\it timelike Szab\'o}, {\it spacelike Ivanov-Petrova}, and {\it timelike
Ivanov-Petrova} are defined similarly. Suppose that $p\ge1$ and $q\ge1$ so the conditions timelike Osserman and spacelike Osserman
are both non-trivial. One can then use analytic continuation to see these two conditions are equivalent. Similarly,
spacelike Szab\'o and timelike Szab\'o are equivalent notions if $p\ge1$ and $q\ge1$. Finally, spacelike Ivanov-Petrova and timelike
Ivanov-Petrova are equivalent notions if
$p\ge2$ and
$q\ge2$. Thus we shall simply speak of {\it Osserman}, {\it Szab\'o}, or {\it Ivanov-Petrova} manifolds; see
\cite{Gi02} for further details.

We shall
refer to \cite{GKV02,Gi02} for a fuller discussion of geometry of the Riemann curvature tensor and shall
content ourselves here with a very brief historical summary. Szab\'o \cite{S91} showed that a Riemannian
manifold is Szab\'o if and only if it is a local symmetric space. Gilkey and Stavrov \cite{GSt02}
showed that a Lorentzian manifold is Szab\'o if and only if it has constant sectional curvature.

Let $\mathcal{M}$ be a Riemannian manifold of dimension $m\ne16$.
Chi \cite{refChi} and Nikolayevsky \cite{N03,N04a,N04b} showed that $\mathcal{M}$ is Osserman if and only if
$\mathcal{M}$ either is  flat or is locally isometric to a rank $1$-symmetric space. This result
settles in the affirmative for
$m\ne16$ a question originally posed by Osserman \cite{O90}. Work of Bla\v zi\'c, Bokan and Gilkey \cite{BBG97} and
of Garc\'{\i}a--R\'{\i}o, Kupeli and  V\'azquez-Abal \cite{GKV97} showed a Lorentzian manifold is Osserman
if and only if it has constant sectional curvature.

Work of of Gilkey \cite{Gi99}, of Gilkey, Leahy, Sadofsky \cite{GLS99}, and of
Nikolayevsky \cite{Ni04c} showed that a Riemannian manifold is Ivanov-Petrova if and only if it either has
constant sectional curvature or it is locally isometric to a warped product of an interval $I$ with a metric
of constant sectional curvature $K$ where the warping function
$f(t)=Kt^2+At+B$ is quadratic and non-vanishing for $t\in I$. This result was extended to the Lorentzian setting for $q\ge11$
by Zhang
\cite{Z02}; results of Stavrov \cite{St03} provide some insight into the higher signature setting.

\subsection{Nilpotency} The picture is very different when $p\ge2$ and $q\ge2$ and the classification of Osserman, Ivanov-Petrov, and
Szab\'o manifolds is far from complete.  The eigenvalue $0$ plays a
distinguished role. We say that
$\mathcal{M}$ is {\it nilpotent Osserman} if
$0$ is the only eigenvalue of $J$ or equivalently if $J(\xi)^m=0$ for any tangent vector $\xi$; the notions
{\it nilpotent Szab\'o} and {\it nilpotent Ivanov-Petrova} are defined similarly.

\subsection{Holonomy} Let $\gamma$ be a smooth curve in a pseudo-Riemannian manifold $\mathcal{M}$. Parallel translation along
$\gamma$ defines a linear isometry $P_\gamma:T_{\gamma(0)}M\rightarrow T_{\gamma(1)}M$. The set of all such automorphisms where
$\gamma(0)=\gamma(1)$ forms a group which is called the {\it holonomy group}; we shall denote this group
by $\mathcal{H}_P(\mathcal{M})$.

\subsection{Generalized plane wave manifolds} Let
$x=(x_1,...,x_m)$ be the usual coordinates on
$\mathbb{R}^m$.  We say
$\mathcal{M}:=(\mathbb{R}^m,g)$ is a {\it generalized plane wave manifold} if
$$\nabla_{\partial_{x_i}}\partial_{x_j}=\textstyle\sum_{k>\max(i,j)}\Gamma_{ij}{}^k(x_1,...,x_{k-1})\partial_{x_k}\,.$$
Let
$\mathcal{T}$ be the nilpotent upper triangular group of all matrices of the form:
$$T=\left(\begin{array}{llllll}
1&*&*&...&*&*\\
0&1&*&...&*&*\\
0&0&1&...&*&*\\
...&...&...&...&...&...\\
0&0&0&...&1&*\\
0&0&0&...&0&1\end{array}\right)\,.$$
\begin{theorem}\label{thm-1.1} Let $\mathcal{M}$ be a generalized plane wave manifold. Then:
\begin{enumerate}
\item $\mathcal{M}$ is complete and strongly geodesically convex.
\item $\nabla_{\partial_{x_{j_1}}}...\nabla_{\partial_{x_{j_\nu}}}R(\partial_{x_{i_1}},\partial_{x_{i_2}})\partial_{x_{i_3}}$\newline
\phantom{.}\qquad$
    =\textstyle\sum_{k>\max(i_1,i_2,i_3,j_1,...j_\nu)}R_{i_1i_2i_3}{}^k{}_{;j_1...j_\nu}
   (x_1,...,x_{k-1})\partial_{x_k}$.
\item $\mathcal{M}$ is nilpotent Osserman, nilpotent Ivanov-Petrova, and nilpotent Szab\'o.
\item $\mathcal{M}$ is Ricci flat and Einstein.
\item $\mathcal{M}$ is VSI.
\item If $\gamma$ is a smooth curve in $\mathbb{R}^m$, then
$P_\gamma\partial_{x_i}=\partial_{x_i}+\textstyle\sum_{j>i}a^j\partial_{x_j}$.
\item  $\mathcal{H}_P(\mathcal{M})\subset\mathcal{T}$.
\end{enumerate}
\end{theorem}

We shall establish Theorem \ref{thm-1.1} in \S\ref{sect2}. Since all the scalar Weyl invariants vanish, one of the central
difficulties in this subject is constructing isometry invariants of such manifolds. In the remaining sections
of this paper, we present several other families of examples with useful geometric properties and exhibit appropriate local
invariants which are not of Weyl type.

\section{Geometric properties of generalized plane wave manifolds}\label{sect2}

\subsection{Geodesics} We begin the proof of Theorem \ref{thm-1.1} by examining the geodesic structure. Let
$\gamma(t)=(x_1(t),...,x_m(t))$ be a curve in
$\mathbb{R}^m$;
$\gamma$ is a geodesic if and only
\begin{eqnarray*}
&&\ddot x_1(t)=0,\ \ \text{and for}\ \ k>1\quad\text{we have}\\
&&\ddot x_k(t)+\textstyle\sum_{i,j<k}\dot x_i(t)\dot
x_j(t)\Gamma_{ij}{}^k(x_1,...,x_{k-1})(t)=0\,.
\end{eqnarray*}
We solve this system of equations recursively. Let
$\gamma(t;\vec x^{\phantom{.}0},\vec x^{\phantom{.}1})$ be defined by
\begin{eqnarray*}
&&x_1(t):=x_1^0+x_1^1t,\ \ \text{and for}\ \ k>1\\
&&x_k(t):=x_k^0+x_k^1t-\textstyle\int_0^t\int_0^s\textstyle\sum_{i,j<k}\dot x_i(r)\dot
x_j(r)\Gamma_{ij}{}^k(x_1,...,x_{k-1})(r)drds\,.
\end{eqnarray*}
Then $\gamma(0;\vec x^{\phantom{.}0},\vec x^{\phantom{.}1})=\vec x^{\phantom{.}0}$ while $\dot\gamma(0;\vec
x^{\phantom{.}0},\vec x^{\phantom{.}1})=\vec x^{\phantom{.}1}$. Thus every geodesic arises in this way so all geodesics
extend for infinite time. Furthermore, given $P,Q\in\mathbb{R}^n$, there is a unique geodesic $\gamma=\gamma_{P,Q}$ so that
$\gamma(0)=P$ and
$\gamma(1)=Q$ where
\begin{eqnarray*}
&&x_1^0=P_1,\ \ x_1^1=Q_1-P_1,\ \ \text{and for}\ \ k>1\quad\text{we have}\\
&&x_k^0=P_k,\ \ x_k^1=Q_k-P_k+\textstyle\int_0^1\int_0^s\textstyle\sum_{i,j<k}\dot
x_i(r)\dot x_j(r)\Gamma_{ij}{}^k(x_1,...,x_{k-1})(r)drds\,.
\end{eqnarray*}
This establishes Assertion (1) of Theorem \ref{thm-1.1}.

\subsection{Curvature} We may expand
\begin{eqnarray*}
R_{ijk}{}^l&=&
     \partial_{x_i}\Gamma_{jk}{}^l(x_1,...,x_{l-1})-\partial_{x_j}\Gamma_{ik}{}^l(x_1,...,x_{l-1})\\
    &+&\Gamma_{in}{}^l(x_1,...,x_{l-1})\Gamma_{jk}{}^n(x_1,...,x_{n-1})\\
    &-&\Gamma_{jn}{}^l(x_1,...,x_{l-1})\Gamma_{ik}{}^n(x_1,...,x_{n-1})\,.
\end{eqnarray*}
As we can restrict the quadratic sums to $n<l$,
$R_{ijk}{}^l=R_{ijk}{}^l(x_1,...,x_{l-1})$.
Suppose $l\le k$. Then $\Gamma_{jk}{}^l=\Gamma_{ik}{}^l=0$.
Furthermore for either of the quadratic terms to be non-zero, there must exist an index $n$ with $k<n$ and
$n<l$. This is not possible if $l\le k$. Thus $R_{ijk}{}^l=0$ if $l\le k$. Suppose $l\le i$. Then
$$\partial_{x_i}\Gamma_{jk}{}^l(x_1,...,x_{l-1})=0\quad
\text{and}\quad\partial_{x_j}\Gamma_{ik}{}^l=\partial_{x_j}0=0\,.
$$
We have $\Gamma_{in}{}^l=0$. For the other quadratic term to be non-zero, there must exist an index $n$ so $i<n$ and $n<l$. This
is not possible if $l\le i$. This shows $R_{ijk}{}^l=0$ if $l\le i$; similarly $R_{ijk}{}^l=0$ if $l\le j$.

This establishes
Assertion (2) of Theorem \ref{thm-1.1} if $\nu=0$, i.e. for the undifferentiated curvature tensor $R$. To study $\nabla R$, we expand
\begin{eqnarray}
R_{ijk}{}^n{}_{;l}&=&\partial_lR_{ijk}{}^n(x_1,...,x_{n-1})\label{eqn-2.a}\\
&-&\textstyle\sum_rR_{rjk}{}^n(x_1,...,x_{n-1})\Gamma_{li}{}^r(x_1,...,x_{r-1})\label{eqn-2.b}\\
&-&\textstyle\sum_rR_{irk}{}^n(x_1,...,x_{n-1})\Gamma_{lj}{}^r(x_1,...,x_{r-1})\label{eqn-2.c}\\
&-&\textstyle\sum_rR_{ijr}{}^n(x_1,...,x_{n-1})\Gamma_{lk}{}^r(x_1,...,x_{r-1})\label{eqn-2.d}\\
&-&\textstyle\sum_rR_{ijk}{}^r(x_1,...,x_{r-1})\Gamma_{lr}{}^n(x_1,...,x_{n-1})\label{eqn-2.e}\,.
\end{eqnarray}
To see $R_{ijk}{}^n{}_{;l}=R_{ijk}{}^n{}_{;l}(x_1,...,x_{n-1})$, we observe that we have:
\begin{enumerate}
\item $i<r<n$ in (\ref{eqn-2.b});
\item $j<r<n$ in (\ref{eqn-2.c});
\item $k<r<n$ in (\ref{eqn-2.d});
\item $r<n$ in (\ref{eqn-2.e}).
\end{enumerate}
To show $R_{ijk}{}^n{}_{;l}=0$ if $n\le\max(i,j,k,l)$, we note that
\begin{enumerate}
\item  $\partial_lR_{ijk}{}^n(x_1,...,x_{n-1})=0$ if $n\le\max(i,j,k,l)$ in (\ref{eqn-2.a});
\item $n>\max(r,j,k)$ and $r>\max(i,l)$ so $n>\max(i,j,k,l)$ in (\ref{eqn-2.b});
\item $n>\max(i,r,k)$ and $r>\max(l,j)$ so $n>\max(i,j,k,l)$ in (\ref{eqn-2.c});
\item $n>\max(i,j,r)$ and $r>\max(k,l)$ so $n>\max(i,j,k,l)$ in (\ref{eqn-2.d});
\item $n>\max(l,r)$ and $r>\max(i,j,k)$ so $n>\max(i,j,k,l)$ in (\ref{eqn-2.e}).
\end{enumerate} This establishes Assertion (2) of
Theorem \ref{thm-1.1} if $\nu=1$ so we are dealing with $\nabla R$. The argument is the same for higher values of $\nu$ and is
therefore omitted.

\subsection{The geometry of the curvature tensor} By Assertion (2) of Theorem \ref{thm-1.1},
\begin{eqnarray*}
&&J(\xi)\partial_{x_i}\subset\operatorname{Span}_{k>i}\{\partial_{x_k}\},\quad
\mathcal{S}(\xi)\partial_{x_i}\subset\operatorname{Span}_{k>i}\{\partial_{x_k}\},\\
&&\mathcal{R}(\pi)\partial_{x_i}\subset\operatorname{Span}_{k>i}\{\partial_{x_k}\}\,.
\end{eqnarray*}
Thus $J$, $\mathcal{R}$, and $\mathcal{S}$ are nilpotent which  proves Assertion (3) of Theorem \ref{thm-1.1}. Furthermore, because
$J(\xi)$ is nilpotent,
$\rho(\xi,\xi)=\operatorname{Tr}(J(\xi))=0$. This implies $\rho=0$ which completes the proof of Assertion (4)
of Theorem \ref{thm-1.1}.

\subsection{Local scalar invariants} Let $\Theta$ be a Weyl monomial which is formed by contracting upper and lower indices in
pairs in the variables
$\{g^{ij},g_{ij},R_{i_1i_2i_3}{}^{i_4}{}_{;j_1...}\}$. The single
upper index in $R$ plays a distinguished role. We choose a
representation for $\Theta$ so the number of $g_{ij}$ variables is
minimal; for example, we can eliminate the $g_{i_3i_4}$ variable
in Equation (\ref{eqn-1.b}) by expressing:
$$|R|^2=g^{i_1j_1}g^{i_2j_2}R_{i_1i_2k}{}^lR_{j_2j_1l}{}^k\,.$$

Suppose there is a
$g_{ij}$ variable in this minimal representation, i.e. that
$$\Theta=g_{ij}R_{u_1u_2u_3}{}^i{}_{;...}R_{v_1v_2v_3}{}^j{}_{;...}...\,.$$
Suppose further that $g^{u_1w_1}$ appears in $\Theta$, i.e. that
$$\Theta=g_{ij}g^{u_1w_1}R_{u_1u_2u_3}{}^i{}_{;...}R_{v_1v_2v_3}{}^j{}_{;...}...\,.$$
We could then raise and lower an index to express
$$\Theta=R^{w_1}{}_{u_2u_3j;...}R_{v_1v_2v_3}{}^j{}_{;...}...
  =R_{ju_3u_2}{}^{w_1}{}_{;...}R_{v_1v_2v_3}{}^j{}_{;...}...$$
which has one less $g_{..}$ variable. This contradicts the assumed minimality. Thus $u_1$ must be contracted against
an upper index; a similar argument shows that $u_2$, $u_3$, $v_1$, $v_2$, and $v_3$ are contracted against an upper
index as well. Consequently
$$\Theta=g_{ij}R_{u_1u_2u_3}{}^i{}_{;...}R_{v_1v_2v_3}{}^j{}_{;...}R_{w_1w_2w_3}{}^{u_1}{}_{;...}...\,.
$$
Suppose $w_1$ is not contracted against an upper index. We then have
\begin{eqnarray*}
\Theta&=&g_{ij}g^{w_1x_1}R_{u_1u_2u_3}{}^i{}_{;...}R_{v_1v_2v_3}{}^j{}_{;...}R_{w_1w_2w_3}{}^{u_1}{}_{;...}...\\
&=&R_{u_1u_2u_3j;...}R_{v_1v_2v_3}{}^j{}_{;...}R^{x_1}{}_{w_2w_3}{}^{u_1}{}_{;...}...\\
&=&g^{u_1y_1}R_{u_1u_2u_3j;...}R_{v_1v_2v_3}{}^j{}_{;...}R^{x_1}{}_{w_2w_3y_1;...}...\\
&=&R{}^{y_1}{}_{u_2u_3j;...}R_{v_1v_2v_3}{}^j{}_{;...}R^{x_1}{}_{w_2w_3y_1;...}\\
&=&R_{ju_3u_2}{}^{y_1}{}_{;...}R_{v_1v_2v_3}{}^j{}_{;...}R^{x_1}{}_{w_2w_3y_1;...}
\end{eqnarray*}
which has one less $g_{ij}$ variable. Thus $w_1$ is contracted against
an upper index so
$$\Theta=g_{ij}R_{u_1u_2u_3}{}^i{}_{;...}R_{v_1v_2v_3}{}^j{}_{;...}R_{w_1w_2w_3}{}^{u_1}{}_{;...}
R_{x_1x_2x_3}{}^{w_1}{}_{;...}...\,.$$
We continue in this fashion to build a monomial of infinite length. This is not possible. Thus we can always
find a representation for $\Theta$ which contains no $g_{ij}$ variables in the summation.

We suppose the evaluation of $\Theta$ is non-zero and argue for a
contradiction. To simplify the notation, group all the lower
indices together. By considering the pairing of upper and lower
indices, we see that we can expand $\Theta$ in cycles:
$$\Theta=R_{...i_r...}{}^{i_1}R_{...i_1...}{}^{i_2}...R_{...i_{r-1}...}{}^{i_r}...\,.$$
By Theorem \ref{thm-1.1} (2), $R_{...j...}{}^l=0$ if $l\le j$.
Thus the sum runs over indices where
 $i_r<i_1<i_2<...<i_r$. As this is the empty sum, we see that $\Theta=0$ as desired.
\subsection{Holonomy} Let $X=\sum_ia_i(t)\partial_{x_i}$ be a vector field which is defined along a curve
$\gamma=(\gamma_1,...,\gamma_m)$ in
$\mathbb{R}^m$. Then
$\nabla_{\dot\gamma}X=0$ if and only if
$$
0=\textstyle\sum_i\dot a_i(t)\partial_{x_i}+\textstyle\sum_{i,j,k:i,j<k}\Gamma_{ij}{}^k(t)a_i(t)\dot\gamma_j(t)\partial_{x_k}\,.
$$
Consequently, we can solve these equations by taking recursively
$$
a_k(t)=a_k(0)-\textstyle\int_0^t\textstyle\sum_{i,j<k}\Gamma_{ij}{}^k(a_1(s),...,a_{k-1}(s))a_i(s)\dot\gamma_j(s)ds\,.
$$
If $a_i(0)=0$ for $i<\ell$, we may conclude $a_i(t)=0$ for all $t$ if $i<\ell$. Assertions (6) and (7) now follow. This
completes the proof of Theorem \ref{thm-1.1}.\hfill\qedbox

\section{Manifolds of signature $(2,2+k)$}\label{sect3}

\subsection{The manifolds $\mathcal{M}_{4+k,F}^0$} Let $(x,y,z_1,...,z_k,\tilde y,\tilde x)$ be coordinates on $\mathbb{R}^{4+k}$. Let
$F(y,z_1,...,z_k)$ be an affine function of $(z_1,...,z_k)$, i.e.
$$F(y,z_1,...,z_k)=f_0(y)+f_1(y)z_1+...+f_k(y)z_k\,.$$
Let $\mathcal{M}_{4+k,F}^0:=(\mathbb{R}^{4+k},g^0_{4+k,F})$ where:
\begin{eqnarray*}
&&g^0_{4+k,F}(\partial_x,\partial_{\tilde x})=g^0_{4+k,F}(\partial_y,\partial_{\tilde y})=g^0_{4+k,F}(\partial_{z_i},\partial_{z_i})=1,\\
&&g^0_{4+k,F}(\partial_x,\partial_x)=-2F(y,z_1,...,z_k)\,.
\end{eqnarray*}

\begin{theorem}\label{thm-3.1}
$\mathcal{M}^0_{4+k,F}$ is a generalized plane wave manifold of signature $(2,2+k)$.
\end{theorem}

\begin{proof}The non-zero Christoffel symbols of the first kind are given by
\begin{eqnarray*}
&&g^0_{4+k,F}(\nabla_{\partial_x}\partial_x,\partial_y)=f_0^\prime+\textstyle\sum_if_i^\prime z_i,\\
&&g^0_{4+k,F}(\nabla_{\partial_y}\partial_x,\partial_x)=g^0_{4+k,F}(\nabla_{\partial_x}\partial_y,\partial_x)
   =-\{f_0^\prime+\textstyle\sum_if_i^\prime z_i\},\\
&&g^0_{4+k,F}(\nabla_{\partial_x}\partial_x,\partial_{z_i})=f_i,\\
&&g^0_{4+k,F}(\nabla_{\partial_{z_i}}\partial_x,\partial_x)
   =g^0_{4+k,F}(\nabla_{\partial_x}\partial_{z_i},\partial_x)=-f_i\,.
\end{eqnarray*}
Consequently the non-zero Christoffel symbols of the second kind are given by
\begin{eqnarray*}
&&\nabla_{\partial_x}\partial_x=\{f_0^\prime+\textstyle\sum_if_i^\prime z_i\}\partial_{\tilde y}+
     \textstyle\sum_if_i\partial_{z_i},\\
&&\nabla_{\partial_y}\partial_x=\nabla_{\partial_x}\partial_y=
     -\{f_0^\prime+\textstyle\sum_if_i^\prime z_i\}\partial_{\tilde x},\\
&&\nabla_{\partial_{z_i}}\partial_x=\nabla_{\partial_x}\partial_{z_i}=-f_i\partial_{\tilde x}\,.
\end{eqnarray*}
This has the required triangular form.
\end{proof}

\subsection{$k$-Curvature homogeneity} Let $\mathcal{M}:=(M,g)$ be a pseudo-Riemannian manifold. If $P\in M$, let
$g_P\in\otimes^2T^*_PM$ be the restriction of $g$ to the tangent space
$T_PM$. We use the metric to lower indices and regard $\nabla^kR\in\otimes^{4+k}T^*M$; let $\nabla^kR_P$ be the restriction of
$\nabla^kR$ to
$T_PM$ and let
$$\mathcal{U}^k(\mathcal{M},P):=(T_PM,g_P,R_P,...,\nabla^kR_P)\,.$$
This is a purely algebraic object. Following
Kowalski, Tricerri, and Vanhecke
\cite{KTV91,KTV92}, we say that
$\mathcal{M}$ is {\it
$k$-curvature homogeneous} if given any two points
$P$ and $Q$ of $M$, there is a isomorphism $\Psi_{P,Q}$ from $\mathcal{U}^k(\mathcal{M},P)$ to $\mathcal{U}^k(\mathcal{M},Q)$, i.e.
a linear isomorphism $\Psi_{P,Q}$ from $T_PM$ to $T_QM$ such that
$$
  \Psi_{P,Q}^*g_Q=g_P\quad\text{and}\quad
  \Psi_{P,Q}^*\nabla^iR_Q=\nabla^iR_P\text{ for }0\le i\le k\,.
$$
Similarly, $\mathcal{M}$ is said to be {\it locally
homogeneous} if given any two points $P$ and $Q$, there are neighborhoods $U_P$ and $U_Q$ of $P$ and $Q$,
respectively, and an isometry
$\psi_{P,Q}:U_P\rightarrow U_Q$ such that $\psi_{P,Q}P=Q$. Taking $\Psi_{P,Q}:=(\psi_{P,Q})_*$ shows that locally
homogeneous manifolds are $k$-curvature homogeneous for any $k$.

More generally, we can consider a {\it $k$-model}
$\mathcal{U}^k:=(V,h,A^0,...,A^k)$ where
$V$ is an
$m$-dimensional real vector space, where $h$ is a non-degenerate inner product of signature $(p,q)$ on $V$, and where
$A^i\in\otimes^{4+i}V^*$ has the appropriate universal curvature symmetries. For example, we assume that:
\begin{equation}\label{eqn-3.a}
\begin{array}{l}
A^0(\xi_1,\xi_2,\xi_3,\xi_4)=A^0(\xi_3,\xi_4,\xi_1,\xi_2)=-A^0(\xi_2,\xi_1,\xi_3,\xi_4)\quad\text{and}\\
A^0(\xi_1,\xi_2,\xi_3,\xi_4)+A^0(\xi_2,\xi_3,\xi_1,\xi_4)+A^0(\xi_3,\xi_1,\xi_2,\xi_4)=0\,.
\vphantom{\vrule height 11pt}
\end{array}\end{equation}
We say that {\it $\mathcal{U}^k$ is a $k$-model
for
$\mathcal{M}$} if given any point
$P\in M$, there is an isomorphism $\Psi_P$ from $\mathcal{U}^k(\mathcal{M},P)$ to $\mathcal{U}^k$. Clearly $\mathcal{M}$ is
$k$-curvature homogeneous if and only if $\mathcal{M}$ admits a $k$-model; one may take as the $k$ model
$\mathcal{U}^k:=\mathcal{U}^k(\mathcal{M},P)$ for any $P\in M$.

\subsection{The manifolds $\mathcal{M}_{6,f}^1$} We specialize the construction given above by taking $F=yz_1+f(y)z_2$. Let
$\mathcal{M}_{6,f}^1:=(\mathbb{R}^6,g^1_{6,f})$ where
\begin{equation}\label{eqn-3.b}
\begin{array}{l}
g^1_{6,f}(\partial_x,\partial_{\tilde x})=g^1_{6,f}(\partial_y,\partial_{\tilde y})
  =g^1_{6,f}(\partial_{z_1},\partial_{z_1})=g^1_{6,f}(\partial_{z_2},\partial_{z_2})=1,\quad\text{and}\\
g^1_{6,f}(\partial_x,\partial_x)=-2(yz_1+f(y)z_2)\,.\vphantom{\vrule
height 12pt}
\end{array}\end{equation}
\subsection{An invariant which is not of Weyl type} Set
\begin{equation}
\label{eqn-3.c}
\alpha^1_6(f,P)=\frac{|f^\prime(P)|}{\sqrt{1+(f^\prime(P))^2}}\,.
\end{equation}

\begin{theorem}\label{thm-3.2}
Assume that $f^{\prime\prime}>0$. Then
\begin{enumerate}
\item $\mathcal{M}^1_{6,f}$ is a $0$-curvature homogeneous generalized plane wave manifold.
\item If $\mathcal{U}^1(\mathcal{M}^1_{6,f_1},P_1)$ and
$\mathcal{U}^1(\mathcal{M}^1_{6,f_2},P_2)$ are isomorphic, then\newline
$\alpha^1_6(f_1,P_1)=\alpha^1_6(f_2,P_2)$.
\item $\alpha^1_6\vphantom{\vrule height 10pt}$ is an isometry invariant of this family which is not of Weyl type.
\item $\mathcal{M}^1_{6,f}$ is not $1$-curvature homogeneous.
\end{enumerate}
\end{theorem}

\begin{proof} We use Theorem \ref{thm-3.1} to see that $\mathcal{M}^1_{6,f}$ is a generalized plane wave manifold. Furthermore, up to the
usual
$\mathbb{Z}_2$ symmetries, the computations performed in the proof of Theorem
\ref{thm-3.1} show that the non-zero entries in the curvature tensor are:
\begin{eqnarray*}
&&R(\partial_x,\partial_y,\partial_y,\partial_x)=f^{\prime\prime}z_2,\quad
  R(\partial_x,\partial_y,\partial_{z_1},\partial_x)=1,\quad
  R(\partial_x,\partial_y,\partial_{z_2},\partial_x)=f^\prime\,.
\end{eqnarray*}
We set
\begin{eqnarray*}
&&X:=c_1\{\partial_x-\ffrac12g_{6,f}^1(\partial_x,\partial_x)\partial_{\tilde x}\},\\
&&\tilde X:=c_1^{-1}\partial_{\tilde x},\\
&&Y:=c_2\{\partial_y-\varepsilon_1\partial_{z_1}-\varepsilon_2\partial_{z_2}
         -\ffrac12(\varepsilon_1^2+\varepsilon_2^2)\partial_{\tilde y}\},\\
&&\tilde Y:=c_2^{-1}\partial_{\tilde y},\\
&&Z_1:=c_3\{\partial_{z_1}+f^\prime\partial_{z_2}+(\varepsilon_1+f^\prime\varepsilon_2)\partial_{\tilde y}\},\\
&&Z_2:=c_3\{\partial_{z_2}-f^\prime\partial_{z_1}+(\varepsilon_2-f^\prime\varepsilon_1)\partial_{\tilde y}\}\,.
\end{eqnarray*}
Since $R(\partial_x,\partial_y,\partial_{z_1},\partial_x)=1$ and
$R(\partial_x,\partial_y,\partial_{z_2},\partial_x)\ne0$, we may choose $\varepsilon_1$, $\varepsilon_1$, $c_1$, $c_2$, and
$c_3$ so that
\begin{eqnarray}
&&R(\partial_x,\partial_y,\partial_y,\partial_x)-2\varepsilon_1R(\partial_x,\partial_y,\partial_{z_1},\partial_x)
-2\varepsilon_2R(\partial_x,\partial_y,\partial_{z_2},\partial_x)=0,\label{eqn-3.d}\\
&&R(\partial_x,\partial_y,\partial_y,\partial_x;\partial_y)
   -3\varepsilon_2R(\partial_x,\partial_y,\partial_y,\partial_x;\partial_{z_2})=0,\label{eqn-3.e}\\
&&c_3^2(1+(f^\prime)^2)=1,\label{eqn-3.f}\\
&&c_3(1+(f^\prime)^2)c_1^2c_2=1,\label{eqn-3.g}\\
&&c_3c_1^2c_2^2f^{\prime\prime}=1\,.\label{eqn-3.h}
\end{eqnarray}
We show that $\mathcal{M}^1_{6,f}$ is $0$-curvature homogeneous and complete the proof of Assertion (1) by noting that the possibly
non-zero entries in these tensors are given by:
$$\begin{array}{ll}g^1_{6,f}(X,\tilde X)=g^1_{6,f}(Y,\tilde Y)=1.\\
g^1_{6,f}(Z_1,Z_1)=g^1_{6,f}(Z_2,Z_2)=1&\text{[see equation (\ref{eqn-3.f})]},\vphantom{\vrule height 11pt}\\
R(X,Y,Y,X)=0&\text{[see equation (\ref{eqn-3.d})]},\vphantom{\vrule height 11pt}\\
R(X,Y,Z_1,X)=1&\text{[see equation (\ref{eqn-3.g})]},\vphantom{\vrule height 11pt}\\
R(X,Y,Z_2,X)=0\,.\vphantom{\vrule height 11pt}
\end{array}$$

The possibly non-zero components of $\nabla R$ are:
\begin{eqnarray*}
&&\nabla R(\partial_x,\partial_y,\partial_y,\partial_x;\partial_{z_2})=
   \nabla R(\partial_x,\partial_y,\partial_{z_2},\partial_x;\partial_y)=f^{\prime\prime}>0,\\
&&\nabla R(\partial_x,\partial_y,\partial_y,\partial_x;\partial_y)=f^{\prime\prime\prime}z_2\,.
\end{eqnarray*}
The possibly non-zero components of $\nabla R$ with respect to this basis are given by:
$$\begin{array}{ll}
\nabla R(X,Y,Y,X;Z_1)=\nabla R(X,Y,Z_1,X;Y)=f^\prime&\text{[see equation (\ref{eqn-3.h})]},\vphantom{\vrule height 11pt}\\
\nabla R(X,Y,Y,X;Y)=0&\text{[see equation (\ref{eqn-3.e})]},\vphantom{\vrule height 11pt}\\
\nabla R(X,Y,Y,X;Z_2)=\nabla R(X,Y,Z_2,X;Y)=1&\text{[see equation (\ref{eqn-3.h})]}\,.\vphantom{\vrule height 11pt}
\end{array}$$

We shall say that a basis $\mathcal{B}=\{{}^1X,{}^1Y,{}^1Z_1,{}^1Z_2,{}^1\tilde Y,{}^1\tilde X\}$ is {\it normalized} if the non-zero
entries in $R$ and
$\nabla R$ are
\begin{eqnarray*}
&&R({}^1X,{}^1Y,{}^1Z_1,{}^1X)=1,\quad\text{and}\\
&&\nabla R({}^1X,{}^1Y,{}^1Y,{}^1X;{}^1Z_2)=\nabla R({}^1X,{}^1Y,{}^1Z_2,{}^1X;{}^1Y)=1\,.\vphantom{\vrule height 12pt}
\end{eqnarray*}
For example, $\mathcal{B}=\{X,Y,Z_1-f^\prime Z_2,Z_2,\tilde Y,\tilde X\}$ is a normalized basis. Let
\begin{eqnarray*}
&&\ker(R):=\{\eta:R(\xi_1,\xi_2,\xi_3,\eta)=0\quad\forall \xi_i\},\\
&&\ker(\nabla R):=\{\eta:\nabla R(\xi_1,\xi_2,\xi_3,\xi_4;\eta)=0\text{ and }
  \nabla R(\xi_1,\xi_2,\xi_3,\eta;\xi_4)=0\quad\forall\xi_i\}\,.
\end{eqnarray*}
It is then immediate that
$$\ker(R)=\operatorname{Span}\{Z_2,\tilde X,\tilde Y\}\quad\text{and}\quad
\ker(\nabla R)=\operatorname{Span}\{Z_1-f^\prime Z_2,\tilde X,\tilde Y\}\,.
$$
Let
$\mathcal{B}:=\{{}^1X,{}^1Y,{}^1Z_1,{}^1Z_2,{}^1\tilde Y,{}^1\tilde X\}$ be any normalized basis. Since
${}^1Z_1\in\ker(\nabla R)$ and ${}^1Z_2\in\ker(R)$, we may expand:
\begin{eqnarray*}
&&{}^1Z_1=a_1(Z_1-f^\prime Z_2)+a_2\tilde X+a_3\tilde Y,\\
&&{}^1Z_2=b_1Z_2+b_2\tilde X+b_3\tilde Y\,.
\end{eqnarray*}
Thus we may compute
$$\frac{|g^1_{6,f}({}^1Z_1,{}^1Z_2)|}{|{}^1Z_1|\cdot|{}^1Z_2|}(P)=\frac{|f^\prime|}{\sqrt{1+(f^\prime)^2}}(P)=\alpha^1_6(f,P)\,.$$
This shows $\alpha^1_6(f,P)$ is an invariant of the $1$-model and establishes Assertion (2).

If $\mathcal{M}^1_{6,f}$ is curvature $1$-homogeneous, then necessarily $\alpha^1_6(f)$ is constant or, equivalently,
$(f^\prime)^2=c(1+(f^{\prime})^2)$ for some constant $c$. Since $(f^\prime)^2<(1+(f^\prime)^2)$, $c<1$. Thus we can
solve for
$(f^{\prime})^2$ to see
$(f^{\prime})^2=\frac c{1-c}$ is constant. This contradicts the assumption $f^{\prime\prime}\ne0$.
\end{proof}

\subsection{Weak curvature homogeneity} We can weaken the notion of curvature homogeneity slightly. Let $A^0\in\otimes^4V^*$ be an
algebraic curvature tensor, i.e.
$A^0$ has the usual symmetries of the curvature tensor given in Equation (\ref{eqn-3.a}).
We say that $\mathcal{M}^1$ is {\it weakly $0$-curvature homogeneous} if for every point $P\in M$, there is an isomorphism
$\Phi:T_PM\rightarrow V$ so that $\Phi^*A^0=R$. There is no requirement that $\Phi$ preserve an inner product. The notion of {weakly
$k$-curvature homogeneous} is similar; we consider models $(V,A^0,...,A^k)$ where $A^i\in\otimes^{4+i}(V^*)$ has the appropriate
curvature symmetries. Since we have lowered all the indices, this is a different notion from the notion of {\it affine $k$-curvature
homogeneity} that will be discussed presently.

The following is an immediate consequence of the arguments given above:

\begin{corollary}\label{cor-3.3}
The manifold $\mathcal{M}^1_{6,f}$ is weakly $1$-curvature homogeneous but not $1$-curvature homogeneous.
\end{corollary}

\subsection{Affine geometry}\label{sect-4.7}
Let $\nabla$ be a torsion free connection on $TM$. Since we do not have a metric, we can
not raise and lower indices. Thus we must regard $\nabla^i$ as a $(i+2,1)$ tensor; instead of working with the tensor
$R_{i_1i_2i_3i_4;j_1...}$, we work with
$R_{i_1i_2i_3}{}^{i_4}{}_{;j_1...}$. We say that
$(M,\nabla)$ is affine {\it $k$-curvature homogeneous} if given any two points $P$ and $Q$ of $M$, there is a linear isomorphism
$\phi:T_PM\rightarrow T_QM$ so that
$\phi^*\nabla^iR_Q=\nabla^iR_P$ for $0\le i\le k$. Taking $\nabla$ to be the
Levi-Civita connection of a pseudo-Riemannian metric then yields that any $k$-curvature homogeneous manifold is necessarily
affine $k$-curvature homogeneous by simply forgetting the requirement that $\phi$ be an isometry; there is no metric present in the
affine setting.
We refer to Opozda \cite{O96,O97} for a further discussion of the subject.  The relevant models are:
\begin{eqnarray*}
&&\mathcal{A}^k(\mathcal{M},P):=(T_PM,R_P,\nabla R_P,...,\nabla^kR_P),\quad\text{where}\\
&&\nabla^iR_P\in\otimes^{3+i}T_PM^*\otimes T_PM\,.
\end{eqnarray*}

In fact the invariant $\alpha_6^1$ is an {\it affine invariant}.
We use note that:
$$\begin{array}{ll}
R(X,Y)Z_1=\tilde X,&R(X,Y)X=-Z_1,\\
R(X,Z_1)Y=\tilde X,&R(X,Z_1)X=-\tilde Y,\\
\nabla_{Z_1}R(X,Y)Y=f^\prime\tilde X,&\nabla_{Z_2}R(X,Y)Y=\tilde X,\\
\nabla_{Z_1}R(X,Y)X=-f^\prime\tilde Y,&\nabla_{Z_2}R(X,Y)X=-\tilde Y,\\
\nabla_YR(X,Y)Z_1=f^\prime\tilde X,&\nabla_YR(X,Y)Z_2=\tilde X,\\
\nabla_YR(X,Z_1)Y=f^\prime\tilde X,&\nabla_YR(X,Z_2)Y=\tilde X,\\
\nabla_YR(X,Z_1)X=-\tilde  Y,&\nabla_YR(X,Z_2)X=-\tilde  Y,\\
\nabla_YR(X,Y)X=-f^\prime Z_1-Z_2\,.\\
\end{array}$$
We define the following subspaces:
\begin{eqnarray*}
W_1:&=&\operatorname{Range}(R)=\operatorname{Span}\{R(\xi_1,\xi_2)\xi_3:\xi_i\in\mathbb{R}^6\},\\
W_2:&=&\operatorname{Range}(\nabla R)=\operatorname{Span}\{\nabla_{\xi_1}R(\xi_2,\xi_3)\xi_4:\xi_i\in\mathbb{R}^6\},\\
W_3:&=&\operatorname{Span}\{R(\xi_1,R(\xi_2,\xi_3)\xi_4)\xi_5:\xi_i\in\mathbb{R}^6\},\\
W_4:&=&\ker(R)=\{\eta\in\mathbb{R}^6:R(\xi_1,\xi_2)\eta=0\ \forall\ \xi_i\in\mathbb{R}^6\},\\
W_5:&=&\ker(\nabla R)=\{\eta\in\mathbb{R}^6:\nabla_{\xi_1}R(\xi_2,\xi_3)\eta=0\ \forall\ \xi_i\in\mathbb{R}^6\}\,.
\end{eqnarray*}
\begin{lemma}\label{lem-3.4} We have
\begin{enumerate}
\item $W_1=\operatorname{Span}\{\tilde X,\tilde Y ,Z_1\}$,
\item $W_2=\operatorname{Span}\{\tilde X,\tilde Y,f^\prime Z_1+Z_2\}$,
\item $W_3=\operatorname{Span}\{\tilde X,\tilde Y\}$,
\item $W_4=\operatorname{Span}\{\tilde X,\tilde Y,Z_2\}$,
\item $W_5=\operatorname{Span}\{\tilde X,\tilde  Y,Z_1-f^\prime Z_2\}$.
\item If $\mathcal{A}^1(\mathcal{M}_{6,f_1}^6,P_1)$ and $\mathcal{A}^1(\mathcal{M}_{6,f_2}^6,P_2)$ are isomorphic,
then\newline $\alpha_6^1(f_1,P_1)=\alpha_6^1(f_2,P_2)$.
\end{enumerate}\end{lemma}

\begin{proof} Assertions (1) and (2) are immediate. We compute
\begin{eqnarray*}
&&R(X,R(X,Y)X)X=R(X,-Z_1)X=\phantom{-}\tilde Y,\\
&&R(X,R(X,Y)X)Y=R(X,-Z_1)Y=-\tilde X,\quad\text{ so }\quad
  \operatorname{Span}\{\tilde X,\tilde Y\}\subset W_3\,.
\end{eqnarray*}
We establish Assertion (3) by establishing the reverse inclusion:
$$
R(\xi_1,R(\xi_2,\xi_3)\xi_4)\xi_5=R(\xi_1,aZ_1+b\tilde X+c\tilde Y)\xi_5
   =R(dX,aZ_1)\xi_5\in\operatorname{Span}\{\tilde X,\tilde Y\}\,.
$$

It is clear $W_4\subset\operatorname{Span}\{\tilde X,\tilde Y,Z_2\}$. Let $\eta=aX+bY+cZ_1+dZ_2+e\tilde X+f\tilde Y\in W_4\,.
$
As $R(X,Y)\eta=0$, we have $-aZ_1+c\tilde X=0$ so $a=0$ and $c=0$. As $R(X,Z_1)\eta=0$, we have $-a\tilde Y+b\tilde X=0$
so $b=0$ as well. Assertion (4) now follows.

It is clear $W_5\subset\operatorname{Span}\{\tilde X,\tilde Y,Z_1-f^\prime Z_2\}$. Let $\eta$ be as above. As
$\nabla_{Z_2}R(X,Y)\eta=0$, $-a\tilde Y+b\tilde X=0$ so $a=b=0$. Since $\nabla_YR(X,Y)\eta=0$,
$(cf^\prime+d)=0$ so $d=-cf^\prime$; this establishes Assertion (5).

Suppose we have an isomorphism from $\mathcal{A}^1(\mathcal{M}_{6,f_1}^6,P_1)$ to $\mathcal{A}^1(\mathcal{M}_{6,f_2}^6,P_2)$. We
ignore the $X$ and $Y$ variables. Then we have an isomorphism $\phi$ from $\mathbb{R}^6$ to itself so that
$\phi(W_i(f_1,P_1))=W_i(f_2,P_2)$ for $1\le i\le 5$. We can work in the spaces $W_i/W_3$ to see that we must have the relations:
\begin{eqnarray*}
&&\phi(Z_1)=a_1Z_1,\quad
\phi(f_1^\prime Z_1+Z_2)=a_2(f_2^\prime Z_1+Z_2),\\
&&\phi(Z_2)=a_3Z_2,\quad
\phi(Z_1-f^\prime Z_2)=a_4(Z_1-f_2^\prime Z_2)\,.
\end{eqnarray*}
This yields $a_1f_1^\prime Z_1+a_3Z_2=a_2f_2^\prime Z_1+a_2Z_2$ and
 $a_1Z_1-a_3f_1^\prime Z_2=a_4Z_1-a_4f_2^\prime Z_2$.
Thus $a_1=a_4$ and $a_3=a_2$ so
$a_1f_1^\prime=a_2f_2^\prime$ and $a_2f_1^\prime=a_1f_2^\prime$. Consequently,
$$a_1a_2f_1^\prime f_1^\prime=a_2a_1f_2^\prime f_2^\prime\,.$$
Since the coefficients $a_i$ are non-zero, the desired conclusion follows.
\end{proof}

\section{Neutral signature generalized plane wave manifolds}\label{sect4}

\subsection{The manifolds $\mathcal{M}^2_{2p,\psi}$} Let $p\ge 2$. Introduce coordinates $(x_1,...,x_p,y_1,...,y_p)$ on
$\mathbb{R}^{2p}$. Let $\psi(x)$ be a symmetric $2$-tensor field on $\mathbb{R}^p$. We define a neutral signature metric $g^2_{2p,\psi}$ on
$\mathbb{R}^{2p}$ and a corresponding pseudo-Riemannian manifold $\mathcal{M}^2_{2p,\psi}$ by:
\begin{eqnarray*}
g^2_{2p,\psi}(\partial_{x_i},\partial_{x_j})=\psi_{ij}(x),\ \
g^2_{2p,\psi}(\partial_{x_i},\partial_{y_j})=\delta_{ij},\ \ \text{and}\ \
g^2_{2p,\psi}(\partial_{y_i},\partial_{y_j})=0\,.
\end{eqnarray*}
\begin{theorem}\label{thm-4.1}
$\mathcal{M}^2_{2p,\psi}$ is a generalized plane wave manifold of signature $(p,p)$.
\end{theorem}
\begin{proof} The non-zero Christoffel symbols of the first kind are given by:
$$
\Gamma_{ijk}^x:=g^2_{2p,\psi}(\nabla_{\partial_{x_i}}\partial_{x_j},\partial_{x_k})=\ffrac12\{\partial_{x_j}\psi_{ik}
+\partial_{x_i}\psi_{jk}-\partial_{x_k}\psi_{ij}\}\,.
$$
>From this, it is immediate that:
$$
\nabla_{\partial_{x_i}}\partial_{x_j}=\textstyle\sum_k\Gamma_{ij}^x{}^k(x)\partial_{y_k}\,.
$$
We set $x_{p+i}=y_i$ to see $\mathcal{M}^2_{2p,\psi}$ is a generalized plane wave
manifold.
\end{proof}

\subsection{Holonomy} The manifolds $\mathcal{M}_{2p,\psi}^2$ present a special case. Let $\oo(p)$ be the Lie algebra of the
orthogonal group; this is the additive group of all skew-symmetric $p\times p$ real matrices. If $A_p$ is such a matrix, let
$\mathcal{G}_{2p}$ be the set of all matrices of the form
$$G(A_p)=\left(\begin{array}{ll}I_p&A_p\\0&I_p\end{array}\right)\,.$$
The map $A_p\rightarrow G(A_p)$ identifies $\oo(p)$ with a subgroup of the upper triangular matrices.

\begin{lemma}\label{lem-4.2}
$\mathcal{H}_{P}(\mathcal{M}_{2p,\psi}^2)\subset\oo(p)$.
\end{lemma}

\begin{proof} Let $\gamma$ be a closed loop in $\mathbb{R}^{2p}$. Let $H_\gamma\partial_{x_i}=X_i$ and
$H_\gamma\partial_{y_i}=Y_i$. Since
$\nabla\partial_{y_i}=0$, $Y_i=\partial_{y_i}$. Expand $X_i=\textstyle\sum_j(a_{ij}\partial_{x_j}+b_{ij}\partial_{y_j})$. Since
$H_\gamma$ is an isometry,
$$g_{2p,\psi}^2(X_i,X_j)=\psi_{ij},\quad
  g_{2p,\psi}^2(X_i,Y_j)=\delta_{ij},\quad\text{and}\quad
  g_{2p,\psi}^2(Y_i,Y_j)=0\,.
$$
The relation $g_{2p,\psi}^2(X_i,Y_j)=\delta_{ij}$ and the observation that $Y_i=\partial_{y_i}$ shows
that $a_{ij}=\delta_{ij}$. Thus
$$g_{2p,\psi}^2(X_i,X_j)=\psi_{ij}+b_{ij}+b_{ji}=\psi_{ij}\,.$$
This shows $b\in\oo(p)$.
\end{proof}

\subsection{Jordan normal form} The eigenvalue structure does not determine the Jordan normal form of a self-adjoint or of a
skew-adjoint endomorphism if the metric is indefinite. We say that
$\mathcal{M}$ is {\it spacelike} (resp. {\it timelike}) {\it Jordan Osserman} if
the Jordan normal form of the Jacobi operator $J$ is constant on the pseudo-sphere bundles of
spacelike (resp. timelike) unit vectors. These two notions are not equivalent. The notions
{\it spacelike Jordan Ivanov-Petrova}, {\it timelike Jordan Ivanov-Petrova}, {\it spacelike Jordan
Szab\'o}, and {\it timelike Jordan Szab\'o} are defined similarly. There are no known examples of spacelike or timelike Jordan
Szab\'o manifolds which are not locally symmetric; $\mathcal{S}(\cdot)$ vanishes identically if and only if $\nabla R=0$.

\subsection{The manifolds $\mathcal{M}^3_{2p,f}$} Let $f(x_1,...,x_p)$
be a smooth function on $\mathbb{R}^p$ and let
$\mathcal{M}^3_{2p,f}:=(\mathbb{R}^{2p},g^3_{2p,f})$ where $g^3_{2p,f}$ is defined by
$\psi_{ij}:=\partial_{x_i}f\cdot\partial_{x_j}f$, i.e.
\begin{eqnarray*}
&&g^3_{2p,f}(\partial_{x_i},\partial_{y_j})=\delta_{ij},\quad
g^3_{2p,f}(\partial_{y_i},\partial_{y_j})=0,\quad\text{and}\\
&&g^3_{2p,f}(\partial_{x_i},\partial_{x_j})=\partial_{x_i}(f)\cdot\partial_{x_j}(f)\,.
\end{eqnarray*}
Let $H_{f,ij}:=\partial_{x_i}\partial_{x_j}f$ be the Hessian.  We use
Theorem
\ref{thm-4.1} and results of Gilkey, Ivanova, and Zhang
\cite{GIZ03} to see that:

\begin{theorem}\label{thm-4.3}
Assume that $H_f$ is non-degenerate. Then\begin{enumerate}
\item $\mathcal{M}^3_{2p,f}$ is a generalized plane wave manifold which is isometric to a hypersurface in a flat space of signature
$(p,p+1)$.
\item $\mathcal{M}^3_{2p,f}$ is spacelike and timelike Jordan Ivanov-Petrova.
\item If $p=2$,
then $\mathcal{M}^3_{2p,f}$ is spacelike and timelike Jordan Osserman.
\item If $p\ge3$ and if
$H_f$ is definite, $\mathcal{M}^3_{2p,f}$ is spacelike and timelike Jordan Osserman.
\item If $p\ge3$ and if $H_f$ is indefinite, $\mathcal{M}^3_{2p,f}$ is neither
 spacelike nor timelike Jordan Osserman.
\item The following conditions are equivalent:
\begin{enumerate}\item $f$ is quadratic.
\item $\nabla R=0$.
\item $\mathcal{M}^3_{2p,f}$ is either spacelike or timelike Jordan Szab\'o.
\end{enumerate}\end{enumerate}\end{theorem}

\subsection{An invariant which is not of Weyl type} If $H_f$ is definite, set
\begin{equation}\label{eqn-4.a}
\alpha^3_{2p}(f,P):=\{H_f^{i_1j_1}H_f^{i_2j_2}H_f^{i_3j_3}H_f^{i_4j_4}H_f^{i_5j_5}R(i_1i_2i_3i_4;i_5)R(j_1j_2j_3j_4;j_5)\}(P)
\end{equation}
where $H_f^{ij}$ denotes the inverse matrix and where we sum over repeated indices. One has the following result
of Dunn and Gilkey
\cite{DG04}:

\begin{theorem}\label{thm-4.4}
Let $p\ge3$.
Assume that the Hessian $H_f$ is definite. Then:
\begin{enumerate}\item $\mathcal{M}^3_{2p,f}$ is
$0$-curvature homogeneous.
\item If $\mathcal{U}(\mathcal{M}^3_{2p,f_1},P_1)$ is isomorphic to $\mathcal{U}(\mathcal{M}^3_{2p,f_2},P_2)$,
then\newline $\alpha^3_{2p}(f_1,P_1)=\alpha^3_{2p}(f_2,P_2)$.
\item $\mathcal{M}^3_{2p,f}$ is not locally homogeneous for generic
$f$.
\end{enumerate}
\end{theorem}

\subsection{The manifolds $\mathcal{M}^4_{4,f}$} Let $(x_1,x_2,y_1,y_2)$ be coordinates on $\mathbb{R}^4$. We consider another
subfamily of the examples considered in Theorem
\ref{thm-4.1}. Let $f=f(x_2)$. Let
$$g^4_{4,f}(\partial_{x_1},\partial_{x_1})=-2f(x_2),\quad g^4_{4,f}(\partial_{x_1},\partial_{y_1})=g^4_{4,f}(\partial_{x_2},\partial_{y_2})=1$$
define $\mathcal{M}^4_{4,f}$. Results of Dunn, Gilkey, and Nik\v cevi\'c \cite{DGS04} show:

\begin{theorem}\label{thm-4.5}
Assume that $f^{(2)}$ and $f^{(3)}$ are never vanishing. The manifold $\mathcal{M}^4_{4,f}$ is
a generalized plane wave manifold of neutral signature
$(2,2)$ which is $1$-curvature homogeneous but not symmetric. The following assertions are equivalent:
\begin{enumerate}\item $f^{(2)}=ae^{\lambda y}$ for some $a,\lambda\in\mathbb{R}-\{0\}$.
\item $\mathcal{M}^4_{4,f}$ is homogeneous.
\item $\mathcal{M}^4_{4,f}$ is $2$-curvature homogeneous.
\end{enumerate}
\end{theorem}

\subsection{An invariant which is not of Weyl type} If $f^{(3)}$ is never vanishing, we set
\begin{equation}\label{eqn-4.b}
\alpha^4_{4,p}(f,P):=\frac{f^{(p+2)}\{f^{(2)}\}^{p-1}}{\{f^{(3)}\}^{-p}}(P)\quad\text{for}\quad p=2,3,...\,.
\end{equation}
In the real analytic context,
these form a complete family of isometry invariants that are not of Weyl type. Again, we refer to Dunn, Gilkey, and Nik\v cevi\'c
\cite{DGS04} for:

\begin{theorem}\label{thm-4.6}
Assume that $f_i$ are real analytic functions on $\mathbb{R}$ and that $f_i^{(2)}$ and $f_i^{(3)}$ are
positive for $i=1,2$. The following assertions are equivalent:
\begin{enumerate}\item There exists an isometry
$\phi:(\mathcal{M}^4_{f_1},P_1)\rightarrow(\mathcal{M}^4_{f_2},P_2)$.
\item We have
$\alpha^4_{4,p}(f_1)(P_1)=\alpha^4_{4,p}(f_2)(P_2)$ for $p\ge2$.
\end{enumerate}
\end{theorem}

\subsection{The manifolds $\mathcal{M}^5_{2p+6,f}$} We consider yet another subfamily of the examples considered in Theorem
\ref{thm-4.1}. Introduce coordinates on
$\mathbb{R}^{2p+6}$ of the form
$(x,y,z_0,...,z_p,\bar x,\bar y,\bar z_0,...,\bar z_p)$. Let $\mathcal{M}^5_{2p+6,f}:=(\mathbb{R}^{2p+6},g^5_{2p+6,f})$ be the
pseudo-Riemannian manifold of signature
$(p+3,p+3)$ where:
\begin{eqnarray*}
&&g^5_{2p+6,f}(\partial_{z_i},\partial_{\bar z_j})=\delta_{ij},\ \
  g^5_{2p+6,f}(\partial_x,\partial_{\bar x})=1,\ \
  g^5_{2p+6,f}(\partial_y,\partial_{\bar y})=1,\\
&&g^5_{2p+6,f}(\partial_x,\partial_x)=-2(f(y)+yz_0+y^2z_1+...+y^{p+1}z_p)\,.
\end{eqnarray*}

\subsection{An invariant which is not of Weyl type} If $f^{(p+4)}>0$,
set
\begin{equation}\label{eqn-4.xz}
\alpha^5_{2p+6,k}(f,P):=\frac{f^{(k+p+3)}\{f^{(p+3)}\}^{k-1}}{\{f^{(p+4)}\}^k}(P)\quad\text{for}\quad k\ge2\,.
\end{equation}
The following result follows from work of
Gilkey and Nik\v cevi\'c \cite{GN04b,GN04c}.

\begin{theorem}\label{thm-4.7} Assume that $f^{(p+3)}>0$ and that $f^{(p+4)}>0$. Then:
\begin{enumerate}
\item  $\mathcal{M}^5_{2p+6,f}$ is a generalized plane wave manifold of signature
$(p+3,p+3)$.
\item $\mathcal{M}^5_{2p+6,f}$ is $p+2$-curvature homogeneous.
\item If $k\ge2$ and if $\mathcal{A}^{k+p+1}(M^5_{2p+6,f_1},P_1)$ and
$\mathcal{A}^{k+p+1}(\mathcal{M}^5_{2p+6,f_2},P_2)$ are isomorphic, then
$\alpha^5_{2p+6,k}(f_1,P_1)=\alpha^5_{2p+6,k}(f_2,P_2)$.
\item
$\alpha^5_{2p+6,k}$ is preserved by any affine diffeomorphism and by any isometry.
\item If $f_i$ are real analytic, if $f_i^{(p+3)}>0$, if $f_i^{(p+4)}>0$, and if for all $k\ge2$ we have that
$\alpha^5_{2p+6,k}(f_1,P_1)=\alpha^5_{2p+6,k}(f_2,P_2)$, then there exists an isometry $\phi$ from $\mathcal{M}^5_{2p+6,f_1}$
to $\mathcal{M}^5_{2p+6,f_2}$ with $f(P_1)=P_2$.
\item The following assertions are equivalent:
\begin{enumerate}\item
$\mathcal{M}^5_{2p+6,f}$ is affine $p+3$-curvature homogeneous.
\item $\alpha^5_{2,p}(f)$ is constant.
\item $f^{(p+3)}=ae^{\lambda y}$ for
$a\ne0$ and
$\lambda\ne0$.
\item $\mathcal{M}^5_{2p+6,f}$ is homogeneous.
\end{enumerate}
\end{enumerate}
\end{theorem}

\section{Generalized plane wave manifolds of signature $(2s,s)$}\label{sect5}

\subsection{The manifolds $\mathcal{M}^6_{3s,F}$} Let $s\ge2$. Introduce coordinates $(\vec u,\vec t,\vec v)$ on $\mathbb{R}^{3s}$
for
$$\vec u:=(u_1,...,u_s),\quad\vec t:=(t_1,...,t_s),\quad\text{and}\quad\vec v:=(v_1,...,v_s)\,.
$$
Let $F=(f_1,...,f_s)$ be a collection of smooth real valued
functions of one variable. Let
$\mathcal{M}^6_{3s,F}=(\mathbb{R}^{3s},g^6_{3s,F})$ be the pseudo-Riemannian manifold of signature $(2s,s)$:
\begin{eqnarray*}
&&g^6_{3s,F}(\partial_{u_i},\partial_{u_i})=-2\{f_1(u_1)+...+f_s(u_s)-u_1t_1-...-u_st_s\},\\
&&g^6_{3s,F}(\partial_{u_i},\partial_{v_i})=g^6_{3s,F}(\partial_{v_i},\partial_{u_i})=1,\quad\text{and}\quad
 g^6_{3s,F}(\partial_{t_i},\partial_{t_i})=-1\,.
\end{eqnarray*}
\subsection{An invariant which is not of Weyl type} Define
\begin{equation}\label{eqn-5.xx}
\alpha^6_{3s}(F,P):=\textstyle\sum_{1\le i\le s}\{f_i^{\prime\prime\prime}(u_i)+4u_i\}^2(P)\,.
\end{equation}
We
refer to Gilkey-Nik\v cevi\'c \cite{GN04a} for the proof of the following result:

\begin{theorem}\label{thm-1.3} Let $s\ge3$. Then
\begin{enumerate}
\item $\mathcal{M}^6_{3s,F}$ is a generalized plane wave manifold of signature
$(2s,s)$.
\item $\mathcal{M}^6_{3s,F}$ is $0$-curvature homogeneous.
\item $\mathcal{M}^6_{3s,F}$ is spacelike Jordan Osserman.
\item $\mathcal{M}^6_{3s,F}$ is spacelike Jordan Ivanov-Petrova of rank $4$.
\item $\mathcal{M}^6_{3s,F}$ is not timelike Jordan Osserman.
\item $\mathcal{M}^6_{3s,F}$ is not timelike Jordan Ivanov-Petrova.
\item If $\mathcal{U}^1(\mathcal{M}^6_{3s,F_1},P_1)$ and
$\mathcal{U}^1(\mathcal{M}^6_{3s,F_2},P_2)$ are isomorphic, then\newline $\alpha^6_{3s}(F_1,P_1)=\alpha^6_{3s}(F_2,P_2)$.
\item $\alpha^6_{3s}$ is an isometry invariant.
\item The following assertions are equivalent:
\begin{enumerate}
\item $f_i^{(3)}(u_i)+4u_i=0$ for $1\le i\le s$.
\item $\mathcal{M}^6_{3s,F}$ is a symmetric space.
\item $\mathcal{M}^6_{3s,F}$ is $1$-curvature homogeneous.
\end{enumerate}\end{enumerate}
\end{theorem}

\section*{Acknowledgements} Research of Peter Gilkey partially supported by the Atlantic Association for Research in the
Mathematical Sciences (Canada) and by the Max Planck Institute in the Mathematical Sciences (Leipzig, Germany).
Research of S. Nik\v cevi\'c partially supported by DAAD (Germany), by the Dierks von Zweck foundation (Germany), and by MM 1646
(Serbia). It is a pleasant task to acknowledge helpful comments from Professors U. Simon and D. Aleeksievski. We also acknowledge the
hospitality of the TU (Berlin) where much of this work was done.


\begin{thebibliography}{AAA}
\bibitem{Ni04c} Y. Nikolayevsky, Riemannian manifolds of dimension 7 whose skew-symmetric curvature operator has constant
eigenvalues; math.DG/0311429.
\end{thebibliography}
\end{document}